\documentclass[a4paper,12pt]{article}

\usepackage{amsmath}
\usepackage{amssymb}
\usepackage{latexsym}

\begin{document}


\begin{center}
\Large{Reeb components of leafwise complex foliations 
\vspace{3pt}
\\ 
and their symmetries I\hspace{-.1em}I}
\vspace{8pt}
\end{center}

\begin{center}
Tomohiro HORIUCHI
\end{center}

\vspace*{0.05cm}


\begin{abstract}
We study the group of leafwise holomorphic smooth automorphisms of Reeb components of leafwise complex foliation which are obtained by a certain Hopf construction. In particular, in the case where the boundary holonomy  is infinitely tangent to the identity, we determine the structure of the group of leafwise holomorphic automorphisms. 
\end{abstract}


\setcounter{section}{-1}


\section{\large{Introduction}}

In this article, we continue the study on the symmetries of Reeb components with codimension one leafwise complex foliations, which we started in our previous paper [HM]. In particular, we proceed to study the symmetries of 5-dimensional Reeb components. Recall that a ($p+1$)-dimensional Reeb component is a compact manifold $R = D^{p} \times S^{1}$ with a (smooth) foliation of codimension one, whose leaves are graphs of smooth functions $f : \textrm{int} D^{p} \to \mathbb{R}$ where $\lim_{z \rightarrow \partial D^{p}} f(z) = + \infty$, and a compact leaf which is the boundary $S^{p-1} \times S^{1}$. Here we identify $R$ with $(D^{p} \times \mathbb{R}) / \mathbb{Z}$. 

In [HM], we studied the group of all leafwise holomorphic smooth automorphisms of 3-dimensional Reeb components with leafwise complex foliations. In particular, in the case where the Reeb component is obtained by the Hopf construction and the holonomy tangent to the identity to the infinite order at the boundary, we determined the structure of the group as follows.
\\
\\
\textbf{Theorem} ([HM]) \,\, Let ($R,\mathcal{F},J$) be a 3-dimensional leafwise complex Reeb component as above and $H=\mathbb{C}^{*} /\langle  \lambda \rangle$ be the boundary elliptic curve. Here $\lambda$ ($|\lambda| > 1$) is a complex number, where $\frac{-1}{2 \pi i} \log \lambda$ is the modulus of $H$. Then the group $\textit{Aut}\,(R,\mathcal{F},J)$ of leafwise holomorphic automorphisms of $(R,\mathcal{F},J)$ is isomorphic to the semi-direct product
\vspace{0.1cm}
$$
\mathcal{Z}_{\varphi,\lambda} \rtimes \left\{ (\mathbb{C}^{\ast} \times Z_{\varphi}) \,/\, (\lambda, \varphi)^{\mathbb{Z}} \right\}
\vspace{0.1cm}
$$
where $\mathcal{Z}_{\varphi,\lambda}$ is the space of solutions to a certain functional equation (see Section 2 for the detail) and $Z_{\varphi}$ is the centralizer of $\varphi$ in $\textit{Diff}^{\,\,\infty}([0,\infty))$.
\\

In the case of complex leaf dimension $2$, the boundary leaf of a Reeb component is called a primary Hopf surface. Kodaira [Ko1] classified the primary Hopf surfaces into five types and gave normal form in each case  (Theorem 1.2). We compute the group of leafwise holomorphic automorphisms for each type of the boundary Hopf surface relying on the normal form and obtain the similar results (Theorem 3.9, 3.10, 3.12 and 3.13 to the above).

In order to describe the automorphism groups, we need some results from [HM] on the space of solutions to certain (system of) functional equations concerning flat functions and an expanding diffeomorphism of the half line which are infinitely tangent to the identity at the origin. These are reviewed in Section 2. 


\section{\large{Preliminaries}}


\subsection{Leafwise complex structure}
\vspace{0.1cm}
 
Let $M$ be a $(2n+q)$-dimensional smooth manifold and $\mathcal{F}$ be a smooth foliation of codimension $q$ on $M$ and $p = 2n$ be the dimension of leaves. We refer general basics for foliation theory to [CC].
\\
\\
\textbf{Definition 1.1} (Leafwise complex structure, cf. [MV]) \,\, $(M,\mathcal{F})$ is said to be equipped with a \textit{leafwise complex structure} if there exits a system of local smooth foliated coordinate charts $(U_{\alpha},\phi_{\alpha})$ where $\phi_{\alpha}:U_{\alpha} \rightarrow V_{\alpha} \subset \mathbb{C}^{n} \times \mathbb{R}^{q} = \{ (z_{1}, \cdots, z_{n},x_{1}, \cdots, x_{q}) \}$ is a smooth diffeomorphism onto an open set $V_{\alpha}$ such that the coordinate change
\[
(w_{1}, \cdots, w_{n}, y_{1}, \cdots, y_{q} ) = \gamma_{\beta \alpha}(z_{1},\cdots,z_{n},x_{1},\cdots,x_{q})
\]
is smooth, $y_{j}$'s depend only on $x_{k}$'s $(j,k = 1,\cdots q)$, and when $x_{k}$'s are fixed $w_{l}$'s are holomorphic in $z_{m}$'s $(l,m = 1, \cdots n)$, where $\gamma_{\beta \alpha}:\phi_{\alpha}(U_{\alpha} \cap U_{\beta}) \rightarrow \phi_{\beta}(U_{\alpha} \cap U_{\beta})$. It is equivalent to that the foliation has complex structures vary smoothly in transverse directions. It is eventually equivalent to that the tangent bundle $\tau \mathcal{F}$ to the foliation is equipped with a smooth integrable almost complex structure $J$. We call $(M,\mathcal{F},J)$ a \textit{leafwise complex foliation}.


\subsection{Hopf surfaces}
\vspace{0.1cm}

Let $W$ be the domain $\mathbb{C}^{2} \,\backslash\, \{ O \}$. A compact complex surface is called a \textit{Hopf surface} if its universal covering is biholomorphic to $W$. Especially, a Hopf surface whose fundamental group is infinite cyclic is called a \textit{primary Hopf surface}. Kodaira classified primary Hopf surfaces in [Ko1, 2].
\\
\\
\textbf{Theorem 1.2} (Kodaira\,\,[Ko1,\,2]) \,\, 1)\,\, Any primary Hopf  surface is a quotient space $W / G^{\mathbb{Z}}$ of $W$ with respect to an infinite cyclic group \,$G^{\mathbb{Z}} = \langle G \rangle$ generated by a complex analytic automorphism $G:W \rightarrow W$ of the form $G(z_{1},z_{2}) = (\lambda z_{1} + \tau z_{2}^{p},\, \mu z_{2})$, where $p$ is a positive integer and $\lambda,~\mu,~\tau$ are complex numbers satisfying $  |\lambda| \geq |\mu| > 1$ and $(\lambda - \mu^{p}) \tau = 0$.
\\
2)\,\, A compact complex surface $S$ is biholomorphic to a primary Hopf surface if and only if it is diffeomorphic to $S^{3} \times S^{1}$.


\subsection{Reeb components by Hopf construction}
\vspace{0.1cm}

\textbf{Construction 1.3} (Hopf construction) \,\, Let $\tilde{R}$ be $\mathbb{C}^{2} \times [0,\infty) \,\backslash\, \{ (O,0) \}$ and $\tilde{\mathcal{F}} = \{ \mathbb{C} \times \{ x \} \,;\, x > 0 \} \sqcup \{ W \times \{ 0 \} \}$ be the foliation on $\tilde{R}$ with the natural complex structure $J_{\textrm{std}}$. Let $T$ be a diffeomorphism of $\tilde{R}$ given by
\[
T (z_{1},z_{2},x) = (G \times \varphi)(z_{1},z_{2},x) = (\lambda z_{1} + \tau z_{2}^{p},\, \mu z_{2},\, \varphi(x))
\]
where $p$ is a positive integer, $\lambda,~\mu,~\tau$ are complex numbers as Theorem 1.2 and $\varphi \in \textit{Diff}^{\,\,\infty}([0,\infty))$ is a diffeomorphism of the half line  satisfying $\varphi(x) - x > 0$ for $x > 0$, namely the origin is an expanding unique fixed point. Then the quotient $R = \tilde{R} / T^{\mathbb{Z}}$ has a foliation $\mathcal{F}$ with leafwise complex  structure induced by $\tilde{\mathcal{F}}$. The boundary $H = W / G^{\mathbb{Z}}$ is a primary Hopf surface, and the holonomy along the boundary leaf coincides with $\varphi$. We call $(R, \mathcal{F}, J)$ the \textit{Reeb component with leafwise complex structure} or the \textit{LC Reeb component}.


\section{\large{Functional equations on flat functions}}

In this section, we review the result on the functional equations which we proved in [HM] in order to determine the automorphism groups of a LC Reeb component.

Let $\varphi \in \textit{Diff}^{\,\,\infty}([0,\infty))$ be a diffeomorphism of the half line which is tangent to the infinite order at $x=0$ and satisfies $\varphi(x) - x > 0$ for $x > 0$. Also we fix a complex number $\lambda$ with $|\lambda| > 1$. Let us consider the following (system of) functional equations on $\beta$, $\beta_{1}$ and $\beta_{2} \in C^{\infty}([0,\infty);\mathbb{C})$ concerning $\varphi$ and $\lambda$. If $\lambda$ is a real number, we can consider the same equations for $\beta_{2} \in C^{\infty}([0,\infty);\mathbb{R})$.
\vspace{5pt}
\\
\quad Equation \,(I)\, :\quad $\beta(\varphi(x))=\lambda \beta(x)$.  
\vspace{5pt}
\\
\quad Equation (I\hspace{-.1em}I) :\quad $\beta_1(\varphi(x)) = \lambda \beta_1(x) + \beta_2(x), \;\;\;$  
$\beta_2(\varphi(x))=\lambda \beta_2(x)$.  
\vspace{5pt} 

First consider these equations on $(0,\infty)$. 
Then, Equation (I) has a lot of solutions 
and if we fix any solution 
$\beta^*(x)\in C^\infty((0,\infty);\mathbb{C})$ 
which never vanishes, \textit{i.e}. $\beta^*(x) \ne 0$ for $x>0$. Then each solution corresponds to a smooth function on $S^1=(0,\infty) / \varphi^{\mathbb{Z}}$ by taking $\beta \to \beta/\beta^*$. This gives a bijective correspondence between the space $\mathcal{Z} = \mathcal{Z}_{\varphi,\lambda}$ of solutions  to (I) on $(0, \infty)$ and $C^\infty(S^1; \mathbb{C})$ as vector space. 

Also take the space $\mathcal{S} = \mathcal{S}_{\varphi, \lambda}$ of solutions to Equation (I\hspace{-.1em}I) on $(0, \infty)$. If we assign $\beta_2$ to  a solution $(\beta_1, \beta_2) \in \mathcal{S}$, we obtain the projection $P_2: \mathcal{S} \to \mathcal Z$. Here the kernel of $P_2$ is nothing but $\mathcal{Z}$. We also see that the projection $P_2$ is surjective because for any $\beta_2\in \mathcal Z$
\vspace{0.1cm}
$$
\beta_1(x)= \frac{1}{\lambda \log \lambda} \beta_2(x) \log \beta^*(x)
\vspace*{0.1cm}
$$
gives a solution $(\beta_1, \beta_2)\in \mathcal S$, 
where for $\log\beta^*(x)$ any smooth branch can be taken. 
Therefore, as a vector space,  
$\mathcal S$ has a structure such that 
$$
0 \to \mathcal{Z} \to \mathcal{S} \to \mathcal{Z} \to 0
$$
is a short exact sequence.

We proved the flatness for $\beta \in \mathcal{Z}$ and $(\beta_{1},\beta_{2}) \in \mathcal{S}$ in [HM]. Key points of the proof are the infinte tangency of $\varphi$ and the formula of Fa\'{a} di Bruno.
\\
\\
\textbf{Theorem 2.1} ([HM]) \,\, 1) Any solution $\beta \in \mathcal{Z}$ extends to $[0,\infty)$ so as to be a smooth function which is flat at $x = 0$, \textit{i.e.}\,\,$k$-th jet satisfies $j^{k} \beta(0) = 0$ for any $k = 0, \, 1, \, 2,\, \cdots$.
\\
\,\,2) The same applies to any solution $(\beta_{1},\beta_{2}) \in \mathcal{S}$.
\\
\\
\textbf{Remark 2.2} \,\, Let us consider the following system of functional equations on $\beta_{1}$ and $\beta_{2} \in C^{\infty}([0,\infty);\mathbb{C})$ concerning $\varphi$, $\lambda$ and a non-zero constant $c$.
\vspace{5pt}
\\
\quad Equation (I\hspace{-.1em}I$c$) :\quad $\beta_1(\varphi(x)) = \lambda \beta_1(x) + c \beta_2(x), \;\;\;$  
$\beta_2(\varphi(x))=\lambda \beta_2(x)$.  
\vspace{5pt}
\\
Let us take the space $\mathcal{S}(c) = \mathcal{S}_{\varphi,\lambda}(c)$ of solutions to Equation (I\hspace{-.1em}I$c$) on $(0,\infty)$. Then, by taking $(\beta_{1},\beta_{2}) \mapsto (\beta_{1},c\beta_{2})$, $\mathcal{S}(c)$ is in one-to-one correspondence with $\mathcal{S}$. In particular, any solution $(\beta_{1},\beta_{2}) \in \mathcal{S}(c)$ extend to $[0,\infty)$ so as to be smooth and flat at $x=0$.


\section{\large{Symmetries of 5-dimensional Reeb component}}

In this section we compute the group of automorphisms of a Reeb component of dimension $5$ which is given by a Hopf construction. 
\\
\\
\textbf{Definition 3.1} \,\, Let $(M,\mathcal{F},J)$ be a smooth foliated manifold with leafwise complex structure. A diffeomorphism $f:M \rightarrow M$ is said to be a \textit{leafwise holomorphic smooth automorphism} if and only if it preserves the foliation and give rise to biholomorphism between leaves. We denote by $\textit{Aut}\,(M,\mathcal{F},J)$, is denoted by $\textit{Aut}\,M$ for short, the group of leafwise holomorphic smooth automorphisms of $(M,\mathcal{F},J)$.
\\
\\
Let $(R,\mathcal{F},J)$ be a 5-dimensional LC Reeb component with holonomy $\varphi$ tangent to the identity to the infinite order at the origin and satisfies $\varphi(x) - x > 0$ for $x > 0$. Any element $f \in \textit{Aut}\,R$ has a lift $\tilde{f} \in \textit{Aut}\,(\tilde{R},\tilde{\mathcal{F}},J_{\textrm{std}}) \, (= \textit{Aut}\,\tilde{R})$ which takes the form
\vspace{0.2cm}
\[
\tilde{f}(z_{1},z_{2},x) = (\xi_{1}(z_{1},z_{2},x),\,\xi_{2}(z_{1},z_{2},x),\,\eta(x))
\vspace*{0.2cm}
\]
in $\mathbb{C}^{2} \times \mathbb{R}_{\geq0}$-coordinate. A lift $\tilde{f}$ should commutes with the covering transformation $T$, because, $T \circ \tilde{f} = \tilde{f} \circ T^{k}$ for some $k \in \mathbb{Z}$ but it is easy to see that $k = 1$ when it is restricted to the boundary. Therefore, an element in $\textit{Aut}\, \tilde{R}$ is a lift of some element in $\textit{Aut}\,R$ if and only if it commutes with $T$. Let $\textit{Aut}\,(\tilde{R};T)$ denote the centralizer of $T$ in $\textit{Aut}\, \tilde{R}$, namely, the group of all such lifts.
\\
\\
\textbf{Proposition 3.2} \,\, $\textit{Aut}\,R$ is naturally isomorphic to $\textit{Aut}\,(\tilde{R};T) / T^{\mathbb{Z}}$.

\subsection{Automorphism groups of Hopf surfaces}
\vspace{0.1cm}

Let $\textit{Aut}\,H$ be the group of holomorphic automorphisms of a primary Hopf surface $H$. For any element $h \in \textit{Aut}\,H$, there is a lift $\tilde{h} \in \textit{Aut}\,\tilde{H}$ of $h$ such that it commutes with $G$ because $\tilde{H} = W$ is the universal covering of $H$. Let $\mathit{Aut}\,(\tilde{H};G)$ denote the centralizer of \,$G$\, in $\textit{Aut}\,\tilde{H}$. Then, $\textit{Aut}\,H$ is naturally isomorphic to $\textit{Aut}\,(\tilde{H};G) / G^{\mathbb{Z}}$. Moreover, by Hartogs's theorem, $\tilde{h}$ is extended to an element in $\textit{Aut}\,(\mathbb{C}^{2},\{O\})$ which consists of all automorphisms of $\mathbb{C}^{2}$ fixing the origin. The automorphism $\tilde{h} \times \textit{id}_{\mathbb{R}_{\geq 0}}$ of $\tilde{R}$ clearly commutes with $T$ and defines an element in $\textit{Aut}\,R$. Consequently, we obtain the following.
\\
\\
\textbf{Proposition 3.3} \,\, The restriction map $r_{H}: \mathit{Aut}\,R \to \mathit{Aut}\,H$ is surjective.
\\
\\
By this proposition, the study of the structure of $\mathit{Aut}\,R$ breaks into two parts, that of the kernel $\mathit{Aut}\,(R,H)$ and the restriction map $r_{H}$.
\\

Namba [Na] determined the centralizer of $G$ in $\textit{Aut}\,(\mathbb{C}^{2},\{O\})$ and the automorphism group of $H$ in the case where $G$ is given by the multiplication of a matrix, namely in the cases 1) - 4) below. The following classification 1) - 5) easily follows from Kodaira's result (Theorem 1.2).
\\
\\
\textbf{Theorem 3.4} (Namba [Na]) \,\, An element $F \in \textit{Aut}\,\tilde{H}$ belongs to $\textit{Aut}\,(\tilde{H};G)$ if and only if it is described in the following form.
\vspace{0.3cm}
\\
\,\,1) If $\tau = 0$ and $\lambda = \mu$ ($G(z_{1},z_{2}) = (\lambda z_{1},\, \lambda z_{2})$),
\[
F(z_{1},z_{2}) = (a z_{1} + b z_{2},\,c z_{1} + d z_{2} ),\,\,\, a,\,b,\,c,\,d \in \mathbb{C},\,\,ad-bc \not = 0.
\]
\,\,2) If $\tau = 0$ and $\lambda = \mu^{p}$ for some $p \geq2$ ($G(z_{1},z_{2})=(\mu^{p} z_{1},\,\mu z_{2})$),
\[
F(z_{1},z_{2}) = (a z_{1} + b z_{2}^{p},\, d z_{2} ),\,\,\, a,\,b,\,d \in \mathbb{C},\,\,ad \not = 0.
\]
\,\,3) If $\tau = 0$ and $\lambda \not = \mu^{m}$ for all $m \in \mathbb{N}$ ($G(z_{1},z_{2}) = (\lambda z_{1},\,\mu z_{2})$),  
\[
F(z_{1},z_{2}) = (a z_{1},\,d z_{2} ),\,\,\, a,\,d \in \mathbb{C},\,\,ad \not = 0.
\]
\,\,4) If $\tau \not = 0$ and $\lambda = \mu$ ($G(z_{1},z_{2}) = (\lambda z_{1} + \tau z_{2},\,\lambda z_{2})$),
\[
F(z_{1},z_{2}) = (a z_{1} + b z_{2},\,a z_{2} ),\,\,\, a,\,b \in \mathbb{C},\,\,a \not = 0.
\]
\,\,5) If $\tau \not = 0$ and $\lambda = \mu^{p}$ for some $p \geq 2$ ($G(z_{1},z_{2}) = (\mu^{p} z_{1} + \tau z_{2}^{p},\,\mu z_{2})$),
\[
F(z_{1},z_{2}) = (a^{p} z_{1} + b z_{2}^{p},\,a z_{2} ),\,\,\, a,\,b \in \mathbb{C},\,\,a \not = 0.
\]
\\
\textbf{Corollary 3.5} \,\,\, The group $\mathit{Aut}\,(\tilde{H};G)$ is isomorphic to the following.
\vspace{0.2cm}
\\
\,\,1) If $\tau = 0$ and $\lambda = \mu$, $\textit{Aut}\,(\tilde{H};G) \cong GL(2,\mathbb{C})$.
\vspace{0.1cm}
\\
\,\,2) If $\tau = 0$ and $\lambda = \mu^{p}$ for some $p \geq2$, $\textit{Aut}\,(\tilde{H};G) \cong \mathbb{C} \rtimes (\mathbb{C}^{*} \times \mathbb{C}^{*})$.
\vspace{0.1cm}
\\
\,\,3) If $\tau = 0$ and $\lambda \not = \mu^{m}$ for all $m \in \mathbb{N}$, $\textit{Aut}\,(\tilde{H};G) \cong \mathbb{C}^{\ast} \times \mathbb{C}^{\ast}\subset GL(2,\mathbb{C})$.
\vspace{0.1cm} 
\\
\,\,4) If $\tau \not = 0$ and $\lambda = \mu$, $\textit{Aut}\,(\tilde{H};G) \cong \mathbb{C} \rtimes \mathbb{C}^{*} \subset GL(2,\mathbb{C})$.
\vspace{0.1cm}
\\
\,\,5) If $\tau \not = 0$ and $\lambda = \mu^{p}$ for some $p \geq 2$, $\textit{Aut}\,(\tilde{H};G) \cong \mathbb{C} \rtimes \mathbb{C}^{*}$.
\\
\\
\textbf{Remark 3.6} \,\,\, $\textit{Aut}\,(\tilde{H};G)$ is a subgroup of $GL(2,\mathbb{C})$ in the case 3) and 4). The case 5) is proven in the same way as the case 4). It is also obtained by computing $\textit{Aut}\,(\tilde{R};T)$ like Case 5 in Section 3.2.2.   

\subsection{Structure of $\mathit{Aut}\,(R,H)$}
\vspace{0.1cm}

Here we determine the kernel $\textit{Aut}\,(R,H)$ of the restriction map $r_{H}$. Let $\mathit{Aut}\,(\tilde{R},\tilde{H};T)$ be the subgroup of $\mathit{Aut}\,(\tilde{R};T)$ which consists of all elements which act trivially on the boundary $\tilde{H}$. Any element $f \in \mathit{Aut}\,(R,H)$ has a unique lift to an element in $\mathit{Aut}\,(\tilde{R},\tilde{H};T)$. Namely, 
\\
\\
\textbf{Proposition 3.7} \,\, $\mathit{Aut}\,(R,H)$ is isomorphic to $\mathit{Aut}\,(\tilde{R},\tilde{H};T)$. 
\\
\\
By this proposition, we should compute the group $\textit{Aut}\,(\tilde{R},\tilde{H};T)$. Let us present an element $g \in \mathit{Aut}\,\tilde{R}$ in the form
\vspace{0.2cm}
$$
g(z_{1},z_{2},x) = (\xi_{1}(z_{1},z_{2},x),\, \xi_{2}(z_{1},z_{2},x),\, \eta(x))
\vspace{0.2cm}
$$
where $\xi_{j}:\tilde{R} \to \mathbb{C}$ is holomorphic in $z_{\alpha}$'s ($\alpha = 1,\,2$) and smooth in $x$ for $j = 1,\,2$, and $\eta \in \mathit{Diff}^{\infty}([0,\infty))$. For any $(z_{1},z_{2},x) \in \tilde{R}$,
\vspace{0.2cm}
\[
g \circ T \left( {}^t \! \left( \begin{array}{c} z_{1} \\ z_{2} \end{array} \right) ,\,x \right) = \left( {}^t \! \left( \begin{array}{c} \xi_{1}(\lambda z_{1} + \tau z_{2}^{p},\, \mu z_{2},\, \varphi(x)) \\ \xi_{2}( \lambda z_{1} + \tau z_{2}^{p},\, \mu z_{2},\, \varphi(x)) \end{array} \right),\,\eta(\varphi(x)) \right),
\vspace*{0.3cm}
\]
\[
T \circ g \left( {}^t \! \left( \begin{array}{c} z_{1} \\ z_{2} \end{array} \right) ,\,x \right) = \left( {}^t \! \left( \begin{array}{c} \lambda \xi_{1}(z_{1},z_{2},x) + \tau \xi_{2}(z_{1},z_{2},x)^{p} \\ \mu \xi_{2}(z_{1},z_{2},x) \end{array} \right),\,\varphi(\eta(x)) \right).
\vspace*{0.4cm}
\]
Hence, the element $g=(\xi_{1},\xi_{2},\eta)$ belongs to $\mathit{Aut}\,(\tilde{R};T)$ if and only if it satisfies the following conditions.
\vspace{0.5cm}

\newpage

\hspace{-0.36cm}
Condition (L) :
\vspace{0.1cm}
$$
\xi_{1}(\lambda z_{1} + \tau z_{2}^{p},\, \mu z_{2},\, \varphi(x)) = \lambda \xi_{1}(z_{1},z_{2},x) + \tau \xi_{2}(z_{1},z_{2},x)^{p} \eqno(\textrm{L1})
$$
\vspace{-0.3cm}
$$
\xi_{2}(\lambda z_{1} + \tau z_{2}^{p},\, \mu z_{2},\, \varphi(x)) = \mu \xi_{2}(z_{1},z_{2},x) \eqno(\textrm{L2})
$$
\\
\quad Condition (T) :
\vspace{0.1cm}
$$
\eta \circ \varphi = \varphi \circ \eta,\, \textrm{namely},\, \eta \in Z_{\varphi} = \textrm{the centralizer of} \, \varphi \, \,\textrm{in}\, \mathit{Diff}^{\infty}([0,\infty))
\vspace*{0.3cm}
$$
Furthermore, $g=(\xi_{1},\xi_{2},\eta)$ belongs to $\mathit{Aut}\,(\tilde{R},\tilde{H};T)$ if and only if the above conditions are satisfied with $\xi_{1}(z_{1},z_{2},0) = z_{1}$ and $\xi_{2}(z_{1},z_{2},0) = z_{2}$. For more details on the centralizer $Z_{\varphi}$, see [HM] and the references therein. By assigning $\eta$ to an element $g = (\xi_{1},\xi_{2},\eta) \in \textit{Aut}\,(\tilde{R},\tilde{H};T)$, we obtain the projection $P:\textit{Aut}\,(\tilde{R},\tilde{H};T) \to Z_{\varphi}$. Then, we obtain a short exact sequence
$$
0 \,\to\, \mathcal{K} \,\to\, \textit{Aut}\,(\tilde{R},\tilde{H};T) \,\overset{P}{\to}\, Z_{\varphi} \,\to\, 1 
$$
where $\mathcal{K}$ is the kernel of $P$. In particular, any element $g \in \mathcal{K}$ is described as
$$
g(z_{1},z_{2},x) = (\xi_{1}(z_{1},z_{2},x),\, \xi_{2}(z_{1},z_{2},x),\, x).
\vspace{0.1cm}
$$

If $x = 0$, for $j = 1,\,2$, $\xi_{j}(\cdot,0)$ is holomorphic function on $\tilde{H} = \mathbb{C}^{2} \,\backslash\, \{ O \}$. Since the origin is a removable singularity, $\xi_{j}(\cdot,0)$ is extended  
the holomorphic function on $\mathbb{C}^{2}$ with $\xi_{j}(O,0) = 0$. Then, $\xi_{j}(\cdot,0)$ is expanded in the power series of $z_{1}$ and $z_{2}$ at the origin:
\[
\xi_{j}(z_{1},z_{2},0) = \sum_{k+l > 0} a_{kl}^{j}(0) z_{1}^{k} z_{2}^{l}~~(j = 1,\,2)
\]
Similarly, we expand $\xi_{j}$ in the power series of $z_{1}$ and $z_{2}$ at the origin:
\vspace{0.1cm}
\[
\xi_{j}(z_{1},z_{2},x) = \sum_{k+l \geq 0} a_{kl}^{j}(x) z_{1}^{k} z_{2}^{l}~~(j = 1,\,2)
\]
where $a_{kl}^{j} \in C^{\infty}([0,\infty);\mathbb{C})$ ($k,\,l = 0,\,1,\,2,\,\cdots$), especially, $a_{00}^{j}(0) = 0$.
\\
\quad \, By Proposition 3.3 and 3.7, in order to determine the group of leafwise holomorphic smooth automorphisms of $R$, we redescribe Condition (L) concerning coefficient functions $a_{kl}^{j}$ for any element $g \in\mathit{Aut}\,\tilde{R}$. If $g$ belongs to $\mathit{Aut}\,(\tilde{R},\tilde{H};T)$, $\xi_{j}$ satisfies with $\xi_{j}(z_{1},z_{2},0) = z_{j}$ for $j = 1,\,2$. Then, we may assume that 
\[
a_{10}^{1}(0) = a_{01}^{2}(0) = 1, \,\, a_{01}^{1}(0) = a_{10}^{2}(0) = 0, \,\, a_{kl}^{j}(0) = 0 \,\,\, (j = 1, \, 2, \, k+l \geq 2).
\] 

\subsubsection{Diagonal case $(\tau = 0)$}
\vspace{0.1cm}

Assume $\tau = 0$. Then, we determine the functions $\xi_{1}$ and $\xi_{2}$ satisfying Condition (L). This case is divided as follows. 
\\
\\
\quad \,  \textbf{Case 1}. \, $\tau=0$ \,and\, $\lambda = \mu$.
\\
\\
\quad \, In this case, Condition (L) is written as
\vspace{0.2cm} 
\[
(\textrm{L1}) \,\,\, \xi_{1}(\lambda z_{1}, \lambda z_{2}, \varphi(x)) = \lambda \xi_{1}(z_{1},z_{2},x), \,\,\, (\textrm{L2}) \,\,\, \xi_{2}(\lambda z_{1}, \lambda z_{2}, \varphi(x)) = \lambda \xi_{2} (z_{1},z_{2},x)
\vspace*{0.2cm}
\]
Equation (L1) is reduced to
\vspace{0.2cm}
\[
\sum_{k+l \geq 0} \lambda^{k+l} a_{kl}^{1}(\varphi(x)) z_{1}^{k} z_{2}^{l} = \sum_{k+l \geq 0} \lambda a_{kl}^{1}(x) z_{1}^{k} z_{2}^{l}.
\vspace{0.2cm}
\]
By comparing coefficient functions of $z_{1}^{k} z_{2}^{l}$\,'s term, (L1) is written as the following.
\\
\\
\quad Equation (L1) :
\vspace{0.1cm}
\[
a_{00}^{1}(\varphi(x)) = \lambda a_{00}^{1}(x), \,\, a_{10}^{1}(\varphi(x)) = a_{10}^{1}(x), \,\, a_{01}^{1}(\varphi(x)) = a_{01}^{1}(x),
\vspace{0.1cm}
\]
\[
\lambda^{k+l-1} a_{kl}^{1}(\varphi(x)) = a_{kl}^{1}(x) \,\,\, (k + l \geq 2)
\vspace{0.2cm}
\]
Similarly, (L2) is written as the following.
\\
\\
\quad Equation (L2) :
\vspace{0.1cm}
\[
a_{00}^{1}(\varphi(x)) = \lambda a_{00}^{2}(x), \,\, a_{10}^{2}(\varphi(x)) = a_{10}^{2}(x), \,\, a_{01}^{2}(\varphi(x)) = a_{01}^{2}(x),
\vspace{0.1cm}
\]
\[
\lambda^{k+l-1} a_{kl}^{2}(\varphi(x)) = a_{kl}^{2}(x) \,\,\, (k + l \geq 2)
\vspace{0.1cm}
\]
\\
\textbf{Proposition 3.8} \,\, Let $\varphi \in \textit{Diff}^{\,\,\infty}([0,\infty))$ be a diffeomorphism of the half line which satisfies $\varphi(x) - x > 0$ for $x > 0$, $a \in C^{\infty}([0,\infty);\mathbb{C})$ be a smooth $\mathbb{C}$-valued function on the half line, and $\nu$ be a complex number with $|\nu| > 1$.
\\
\,\,1)\,\, If $a$ satisfies $a(\varphi(x)) = a(x)$ for all $x \in [0,\infty)$, $a(x) = a(0)$ is constant. 
\\
\,\,2)\,\, If $a$ satisfies $\nu a(\varphi(x)) = a(x)$ for all $x \in [0,\infty)$, $a(x) \equiv 0$. 
\\
\\
\textit{Proof}. \,$1)$ \,\, For any $x \geq 0$,\, $a(x) =\displaystyle  \lim_{n \to \infty} a(\varphi^{-n}(x)) = a(0)$.
\vspace{0.1cm}
\\
\,$2)$ \,\, Similarly, for any $x \geq 0$,\, $| a(x) | = \displaystyle  \lim_{n \to \infty} \frac{1}{|\nu|^{n}} |a(\varphi^{-n}(x))| = 0$. 
\hspace{1.5cm}
$\Box$
\vspace{0.3cm}
\\

By Proposition 3.8, for any $x \in [0,\infty)$,
\[
a_{10}^{1}(x) = 1, \quad a_{01}^{1}(x) = 0, \quad a_{10}^{2}(x) = 0, \quad a_{01}^{2}(x) = 1,
\]
\[
a_{kl}^{j}(x) = 0 \,\,\, (j = 1,\,2,\, k+l \geq 2).
\vspace{0.2cm}
\]
Therefore, $g$ belongs to $\mathit{Aut}\,(\tilde{R},\tilde{H};T)$ iff it takes the following form.
\vspace{0.1cm}
\[
g (z_{1},z_{2},x) = (z_{1} + a_{00}^{1}(x),\, z_{2} + a_{00}^{2}(x),\, \eta(x))
\vspace{0.1cm}
\]
where\, $a_{00}^{1}$, $a_{00}^{2} \in \mathcal{Z}_{\varphi,\lambda}$ \,and\, $\eta \in Z_{\varphi}$.
\\

\textbf{Case 2}. \, $\tau = 0$ \,and\, $\lambda = \mu^{p}$ \,for some\, $p \geq 2$.
\\
\\
\quad \, In this case, Condition (L) is written as
\vspace{0.2cm} 
\[
(\textrm{L1}) \,\,\, \xi_{1}(\mu^{p} z_{1}, \mu z_{2}, \varphi(x)) = \mu^{p} \xi_{1}(z_{1},z_{2},x), \,\,\, (\textrm{L2}) \,\,\, \xi_{2}(\mu^{p} z_{1}, \mu z_{2}, \varphi(x)) = \mu \xi_{2} (z_{1},z_{2},x)
\vspace*{0.2cm}
\]
In the same way as Case 1, above equations are reduced to
\vspace{0.2cm}
\[
\sum_{k+l \geq 0} \mu^{pk+l} a_{kl}^{1}(\varphi(x)) z_{1}^{k} z_{2}^{l} = \sum_{k+l \geq 0} \mu^{p} a_{kl}^{1}(x) z_{1}^{k} z_{2}^{l}
\vspace{0.2cm}
\]
\[
\sum_{k+l \geq 0} \mu^{pk+l} a_{kl}^{2}(\varphi(x)) z_{1}^{k} z_{2}^{l} = \sum_{k+l \geq 0} \mu a_{kl}^{2}(x) z_{1}^{k} z_{2}^{l}
\vspace{0.2cm}
\]
By comparing coefficient functions of $z_{1}^{k} z_{2}^{l}$\,'s term, (L1) and (L2) are written as the following.
\\
\\
\quad Equation (L1) :
\vspace{0.1cm}
\[
a_{10}^{1}(\varphi(x)) = a_{10}^{1}(x), \,\, a_{0l}^{1}(\varphi(x)) = \mu^{p-l} a_{0l}^{1}(x) \,\, (0 \leq l \leq p-1),
\vspace{0.1cm}
\]
\[
a_{0p}^{1}(\varphi(x)) = a_{0p}^{1}(x), \,\, \mu^{l-p} a_{0l}^{1}(\varphi(x)) = a_{0l}^{1}(x) \,\,\, (l \geq p+1),
\vspace{0.2cm}
\]
\[
\mu^{p(k-1)+l} a_{kl}^{1}(\varphi(x)) = a_{kl}^{1}(x) \,\,\, (k + l \geq 2,\,k \geq 1)
\vspace{0.3cm}
\]
\quad Equation (L2) :
\vspace{0.1cm}
\[
a_{00}^{2}(\varphi(x)) = \mu a_{00}^{2}(x), \,\, a_{01}^{2}(\varphi(x)) = a_{01}^{2}(x), \,\,  \mu^{p-1} a_{10}^{2}(\varphi(x)) = a_{10}^{2} (x), 
\vspace{0.1cm}
\]
\[
\mu^{pk+l-1} a_{kl}^{2}(\varphi(x)) = a_{kl}^{2}(x) \,\, (k + l \geq 2)
\vspace{0.1cm}
\]
\\
By Proposition 3.8, for any $x \in [0,\infty)$,
\[
a_{10}^{1}(x) = 1, \,\,\, a_{0l}^{1}(x) = 0 \,\,\, (l \geq p), \,\,\, a_{10}^{2}(x) = 0, \,\,\, a_{01}^{2}(x) = 1,
\]
\[
a_{kl}^{1}(x) = 0 \,\,\, (k+l \geq 2, \, k \geq 1), \,\,\,\, a_{kl}^{2}(x) = 0 \,\,\, (k+l \geq 2)
\vspace{0.3cm}
\]
Therefore, $g$ belongs to $\mathit{Aut}\,(\tilde{R},\tilde{H};T)$ iff it takes the following form.
\vspace{0.1cm}
\[
g (z_{1},z_{2},x) = \left( z_{1} + \sum_{j=0}^{p-1} a_{0j}^{1}(x)z_{2}^{j},\, z_{2} + a_{00}^{2}(x),\, \eta(x) \right)
\vspace{0.1cm}
\]
where $a_{0j}^{1} \in \mathcal{Z}_{\varphi,\mu^{p-j}}$ ($j = 0,\,1,\,\cdots,\,p-1$),\, $a_{00}^{2} \in \mathcal{Z}_{\varphi,\mu}$ \,and\, $\eta \in Z_{\varphi}$.
\\

\textbf{Case 3}. \, $\tau = 0$ \,and\, $\lambda \not = \mu^{m}$ \,for all \,$m \in \mathbb{N}$.
\\
\\
\quad \, In this case, Condition (L) is written as
\vspace{0.2cm} 
\[
(\textrm{L1}) \,\,\, \xi_{1}(\lambda z_{1}, \mu z_{2}, \varphi(x)) = \lambda \xi_{1}(z_{1},z_{2},x), \,\,\, (\textrm{L2}) \,\,\, \xi_{2}(\lambda z_{1}, \mu z_{2}, \varphi(x)) = \mu \xi_{2} (z_{1},z_{2},x)
\vspace*{0.2cm}
\]
In the same way as Case 1, above equations are reduced to
\vspace{0.2cm}
\[
\sum_{k+l \geq 0} \lambda^{k} \mu^{l} a_{kl}^{1}(\varphi(x)) z_{1}^{k} z_{2}^{l} = \sum_{k+l \geq 0} \lambda a_{kl}^{1}(x) z_{1}^{k} z_{2}^{l}
\vspace{0.2cm}
\]
\[
\sum_{k+l \geq 0} \lambda^{k} \mu^{l} a_{kl}^{2}(\varphi(x)) z_{1}^{k} z_{2}^{l} = \sum_{k+l \geq 0} \mu a_{kl}^{2}(x) z_{1}^{k} z_{2}^{l}
\vspace{0.2cm}
\]
By comparing coefficient functions of $z_{1}^{k} z_{2}^{l}$\,'s term, (L1) and (L2) are written as the following.
\\
\\
\quad Equation (L1) :
\vspace{0.1cm}
\[
a_{00}^{1}(\varphi(x)) =\lambda a_{00}^{1}(x), \,\, a_{10}^{1}(\varphi(x)) = a_{10}^{1}(x), \,\, \mu^{l} a_{0l}^{1}(\varphi(x)) = \lambda a_{0l}^{1}(x) \,\, (l \geq 1),\vspace{0.1cm}
\]
\[
\lambda^{k-1} \mu^{l} a_{kl}^{1}(\varphi(x)) = a_{kl}^{1}(x) \,\, (k + l \geq 2,\,k \geq 1)
\vspace{0.4cm}
\]
\quad Equation (L2) :
\[
a_{00}^{2}(\varphi(x)) = \mu a_{00}^{2}(x), \,\, a_{01}^{2}(\varphi(x)) = a_{01}^{2}(x), \,\, \lambda^{k} a_{k0}^{2}(\varphi(x)) = \mu a_{k0}^{2}(x) \,\, (k \geq 1),\vspace{0.1cm}
\]
\[
\lambda^{k} \mu^{l-1} a_{kl}^{2}(\varphi(x)) = a_{kl}^{2}(x) \,\, (k + l \geq 2,\,l \geq 1)
\vspace{0.1cm}
\]
\\
By Proposition 3.8\,\,\,2), for any $x \in [0,\infty)$,
\[
a_{kl}^{1}(x) = 0 \,\, (k + l \geq 2,\,k \geq 1),\,\,\, a_{kl}^{2}(x) = 0 \,\,(k+l \geq 1,\,l \geq 1)
\vspace{0.1cm}
\]
Let us consider functional equations concerning $a_{0l}^{1}$ and $a_{k0} ^{2}$. Since $\lambda$ and $\mu$ satisfy the inequality $|\lambda| \geq |\mu| > 1$, they satisfy $|\lambda^{k}| \geq |\mu|$ for $k \geq 1$. Then, we rewrite the functional equation of $a_{k0}^{2}$ as the following. 
\[
\frac{\lambda^{k}}{\mu} a_{k0}^{2} (\varphi(x)) = a_{k0}^{2}(x) \,\,\, (k \geq 1)
\]
If $|\lambda| > |\mu|$, then $a_{k0}^{2}$ is identically equal to $0$ by Proposition 3.8\,\,\,2) because $|\lambda^{k}/\mu|$ is greater than 1 for $k \geq 1$. If $|\lambda| = |\mu|$ and $k = 1$, the equation \vspace{0.1cm}$|a_{10}^{2}(\varphi(x))| = |a_{10}^{2}(x)|$ holds for $x \geq 0$. Then, we have $|a_{10}^{2}(x)| = \displaystyle \lim_{n \to \infty} |a_{10}^{2}(\varphi^{-n}(x))| = |a_{10}^{2}(0)| = 0$. If $|\lambda| = |\mu|$ and $k \geq 2$, then $a_{k0}^{2}$ is also equal to $0$ in the same way as the case $|\lambda| > |\mu|$. 
\vspace{1pt}

Next, we fix the positive integer $p = [\log |\lambda| / \log |\mu|]$, \vspace{1pt}where $[\,\cdot\,]$ is the greatest integer function. Then, $\lambda$ and $\mu$ satisify $|\lambda| \geq |\mu^{l}|$ for $0 \leq l \leq p$ and $|\lambda| < |\mu^{l}|$ for $l \geq p + 1$. The functional equation of $a_{0l}^{1}$ is rewritten as the following.
\[
a_{0l}^{1}(\varphi(x)) = \frac{\lambda}{\mu^{l}} a_{0l}^{1}(x)\,\,\,(0 \leq l \leq p), \,\,\,\, \frac{\mu^{l}}{\lambda} a_{0l}^{1}(\varphi(x)) = a_{0l}^{1}(x)\,\,\,(l \geq p+1) 
\]
If $l \geq p+1$, then $a_{0l}^{1}$ is identically equal to $0$ by Proposition 3.8\,\,\,2) because $|\mu^{l}/\lambda|$ is greater than 1. If $|\lambda| = |\mu^{p}|$, the equation $|a_{0p}^{1}(\varphi(x))| = |a_{0p}^{1}(x)|$ holds for $x \geq 0$, then we have $|a_{0p}^{1}(x)| = \displaystyle \lim_{n \to \infty}|a_{0p}^{1}(\varphi^{-n}(x))| = |a_{0p}^{1}(0)| = 0$.
\vspace{0.1cm} 

Therefore, $g$ belongs to $\mathit{Aut}\,(\tilde{R},\tilde{H};T)$ iff it takes the following form.
\vspace{0.1cm}
\[
g (z_{1},z_{2},x) = \left( z_{1} + \sum_{j=0}^{p} a_{0j}^{1}(x)z_{2}^{j},\, z_{2} + a_{00}^{2}(x),\, \eta(x) \right)
\vspace{0.1cm}
\]
where\, $a_{0j}^{1} \in \mathcal{Z}_{\varphi,\lambda \mu^{-j}}$ ($j = 0,\,1,\,\cdots,\,p$),\, $a_{00}^{2} \in \mathcal{Z}_{\varphi,\mu}$ \,and\, $\eta \in Z_{\varphi}$.
\vspace{0.1cm}
\\

Consequently, we obtain the following theorems.
\\
\\
\textbf{Theorem 3.9} \,\, If $\tau = 0$ and $\lambda = \mu$, then the group $\textit{Aut}(R,H)$ admits a following sequence of extension
\vspace{0.1cm}
\[
0 \,\to\, \mathcal{Z}_{\varphi,\lambda} \times \mathcal{Z}_{\varphi,\lambda} \,\to\, \mathit{Aut}\,(R,H) \cong (\mathcal{Z_{\varphi,\lambda}} \times \mathcal{Z}_{\varphi,\lambda}) \rtimes Z_{\varphi} \,\to\, Z_{\varphi} \,\to\, 1 
\vspace{0.1cm}
\]
where $\mathcal{Z}_{\varphi,\lambda}$ is an infinite dimensional vector space described in Section 2 which is the set of functions $a_{00}^{1}$ and $a_{00}^{2}$ in Case 1 and $\eta \in Z_{\varphi}$ acts on $(\beta_{1},\beta_{2}) \in \mathcal{Z}_{\varphi,\lambda} \times Z_{\varphi,\lambda}$ by $(\beta_{1}(x),\,\beta_{2}(x)) \mapsto (\beta_{1}(\eta(x)),\,\beta_{2}(\eta(x))$. 
\\
\\
\textbf{Theorem 3.10} \,\, If $\tau = 0$ and $\lambda \not = \mu$, then the group $\textit{Aut}\,(R,H)$ admits a following sequence of extension
\vspace{0.1cm}
\[
0 \,\to\, \mathcal{K} \,\to\, \mathit{Aut}\,(R,H) \,\to\, Z_{\varphi} \to 1 
\vspace{0.1cm}
\]
where $\mathcal{K}$ is isomorphic to the following.
\vspace{0.2cm}
\\
\,\,1) If $\lambda = \mu^{p}$ for some $p \geq2$,\, $\mathcal{K}$ admits a sequence of extension
\vspace{0.1cm}
\[
0 \,\to\, \mathcal{Z}_{\varphi,\mu^{p}} \times \mathcal{Z}_{\varphi,\mu^{p-1}} \times \cdots \times \mathcal{Z}_{\varphi,\mu} \,\to\, \mathcal{K} \,\to\, \mathcal{Z}_{\varphi,\mu} \,\to\, 0.
\vspace{0.1cm}
\]
\,\,2) If $\lambda \not = \mu^{m}$ for all $m \geq 1$,\, $\mathcal{K}$ admits a sequence of extension
\vspace{0.1cm}
\[
0 \,\to\, \mathcal{Z}_{\varphi,\lambda} \times \mathcal{Z}_{\varphi, \lambda \mu^{-1}} \times \cdots \times \mathcal{Z}_{\varphi, \lambda \mu^{-p}} \,\to\, \mathcal{K} \,\to\, \mathcal{Z}_{\varphi,\mu} \,\to\, 0.
\vspace{0.1cm}
\]
\quad where $p = \left[ \frac{\log |\lambda|}{\log |\mu|} \right] \in \mathbb{N}$, where $[\,\cdot\,]$ is the greatest integer function.
\vspace{0.1cm}
\\
\\
\textit{Proof}. \,\, $\mathcal{K}$ is regarded as the space of solutions to Condition (L). If $\tau = 0$ and $\lambda \not = \mu$, $\mathcal{K}$ is the set of functions $(a_{00}^{1},\,\cdots,a_{0\,p-1}^{1},\,a_{00}^{2})$ in Case 2 and 3. By assigning $\beta_{2}$ to an element $(\beta_{1,0},\,\cdots,\,\beta_{1,p-1},\,\beta_{2}) \in \mathcal{K}$, we obtain the projection $P_{2}: \mathcal{K} \to \mathcal{Z}_{\varphi,\mu}$. Then, the kernel of $P_{2}$ is $\mathcal{Z}_{\varphi,\mu^{p}} \times \mathcal{Z}_{\varphi,\mu^{p-1}} \times \cdots \times \mathcal{Z}_{\varphi,\mu}$ in Case 2, and it is $\mathcal{Z}_{\varphi,\lambda} \times \mathcal{Z}_{\varphi,\lambda \mu^{-1}} \times \cdots \times \mathcal{Z}_{\varphi,\lambda \mu^{-p}}$ in Case 3. For further details on the group law of $\mathcal{K}$, refer to Remark 3.11.  \hspace{140pt} 
$\Box$
\\
\\
\textbf{Remark 3.11} \,\, For any elements $g_{1}$ and $g_{2} \in \textit{Aut}\,\tilde{R}$ which are described as
\vspace{0.1cm}
\[
g_{i} (z_{1},z_{2},x) = \left( z_{1} + \sum_{j=0}^{p-1} \beta_{1,j}^{(i)}(x)z_{2}^{j},\, z_{2} + \beta_{2}^{(i)}(x),\, x \right) \,\,\, (i = 1,\,2)
\vspace{0.1cm}
\]
where $\beta_{1,0}^{(i)},\, \cdots,\, \beta_{1,p-1}^{(i)},\,\beta_{2}^{(i)} \in C^{\infty}([0,\infty);\mathbb{C})$,
\vspace{15pt}
\\
\,\,\quad $g_{1} \circ g_{2} (z_{1},z_{2},x) = $
\vspace{0.2cm} 
\[
\left( z_{1} + \sum_{j=0}^{p-1} \left\{ \sum_{k=0}^{p-1-j} \left( \begin{array}{c} j + k \\ j \end{array} \right) \beta_{1,j+k}^{(1)} (x) \cdot \beta_{2}^{(2)}(x)^{k} + \beta_{1,j}^{(2)}(x) \right\} z_{2}^{j},\, z_{2} + \{ \beta_{2}^{(1)}(x) + \beta_{2}^{(2)}(x) \},\,x  \right).
\vspace{0.1cm}
\]
In particular, if $p = 1$,
\[
g_{1} \circ g_{2} (z_{1},z_{2},x) = (z_{1} + \{ \beta_{1,0}^{(1)}(x) + \beta_{1,0}^{(2)}(x) \},\, z_{2} + \{ \beta_{2}^{(1)}(x) + \beta_{2}^{(2)} \},\, x ),
\]
and if $\beta_{2}^{(1)} = \beta_{2}^{(2)} = 0$,
\[
g_{1} \circ g_{2} (z_{1},z_{2},x) = \left( z_{1} + \sum_{j=0}^{p-1} \{ \beta_{1,j}^{(1)}(x) + \beta_{1,j}^{(2)}(x)\} z_{2}^{j},\, z_{2},\, x  \right).
\vspace{0.3cm}
\]

\subsubsection{Nondiagonal case~$(\tau \not= 0)$}
\vspace{0.1cm}

Assume $\tau \not =0$. Then, we determine the functions $\xi_{1}$ and $\xi_{2}$ satisfying Condition (L). This case is divided as follows. 
\\
\\
\quad \,  \textbf{Case 4}. \, $\tau \not= 0$ \,and\, $\lambda = \mu$.
\\
\\
\quad \, Let $S$ be the diffeomorphism of $\tilde{R}$ given by $S(z_{1},z_{2},x) = (z_{1}, \tau z_{2},x)$. Then, for any $(z_{1},z_{2},x) \in \tilde{R}$,
\vspace{0.1cm}
\[
S \circ T \circ S^{-1} (z_{1},z_{2},x) = (\lambda z_{1} + z_{2}, \lambda z_{2}, \varphi(x)).
\vspace{0.1cm}
\]
Therefore, a LC Reeb component $R = \tilde{R} / T^{\mathbb{Z}}$ is leafwise holomorphic foliated diffeomorphic to a LC Reeb component $R^{\prime} = \tilde{R} / (S \circ T \circ S^{-1} )^{\mathbb{Z}}$, so we may assume that $\tau$ is equal to 1.

If $\tau = 1$ and ${\lambda = \mu}$, Condition (L) is written as 
\[
(\textrm{L1}) \quad \xi_{1}(\lambda z_{1} + z_{2},\, \lambda z_{2},\, \varphi(x)) = \lambda \xi_{1}(z_{1},z_{2},x) + \xi_{2}(z_{1},z_{2},x), 
\vspace{0.2cm}
\]
\[(\textrm{L2}) \quad \xi_{2}(\lambda z_{1} + z_{2},\, \lambda z_{2},\, \varphi(x)) = \lambda \xi_{2} (z_{1},z_{2},x).
\vspace*{0.3cm}
\]
First, we determine the function $\xi_{2}$ which satisfies (L2). The functional equation of (L2) is reduced to
\[
\sum_{k + l \geq 0} a_{kl}^{2}(\varphi(x)) (\lambda z_{1} +  z_{2})^{k} (\lambda z_{2})^{l} = \sum_{k + l \geq 0} \lambda a_{kl}^{2}(x) z_{1}^{k} z_{2}^{l}.
\]
Furthermore, we compare both coefficient functions of $z_{1}^{k}z_{2}^{l}$'s term.

On $z_{1}^{0}z_{2}^{0}$'s term ($k = l = 0$), the equation is 
\vspace{0.1cm}
\[
a_{00}^{2}(\varphi(x)) = \lambda a_{00}^{2}(x), \,\,\textit{i.e.}\,\, a_{00}^{2} \in \mathcal{Z}_{\varphi,\lambda}.
\vspace{0.1cm}
\]
On $z_{1}^{1}z_{2}^{0}$'s term ($k = 1$, $l =0$), the equation is 
\vspace{0.1cm}
\[
\lambda a_{10}^{2} (\varphi(x)) = \lambda a_{10}^{2}(x).
\vspace{0.1cm}
\]
Then, by Proposition 3.8\,\,\,1), $a_{10}^{2}(x) = a_{10}^{2}(0) = 0$. On $z_{1}^{0}z_{2}^{1}$'s term ($k = 0$, $l=1$), the equation is 
\vspace{0.1cm}
\[
a_{10}^{2}(x) + \lambda a_{01}^{2}(\varphi(x)) = \lambda a_{01}^{2}(x).
\vspace{0.1cm}
\]
Since $a_{10}^{2}$ is equal to 0, above equation is
\vspace{0.1cm}
\[
a_{01}^{2}(\varphi(x)) = a_{01}^{2}(x).
\vspace{0.1cm}
\]
By Proposition 3.8\,\,\,1), $a_{01}^{2}(x) = a_{01}^{2} (0) = 1$. If $k + l \geq 2$, the functional equation of $z_{1}^{k}z_{2}^{l}$'s term is 
\[
\sum_{j=0}^{l}
\left(
\begin{array}{c}
k+l-j
\\
k
\end{array}
\right)
\lambda^{k+l-j} \lambda^{j} a_{k+l-j\,j}^{2}(\varphi(x)) = \lambda a_{kl}^{2}(x).
\]
If $k \geq 2$ and $l=0$, the equation is 
\vspace{0.1cm}
\[
\lambda^{k}a_{k0}^{2}(\varphi(x)) = \lambda a_{k0}^{2}(x).
\vspace{0.1cm}
\]
By Proposition 3.8\,\,\,2), $a_{k0}^{2}(x) = 0$. If $a_{k+l\,0}^{2},\,\cdots,\,a_{k+2\,l-2}$ and $a_{k+1\,l-1}^{2}$ are identically equal to $0$, then the equation of $z_{1}^{k}z_{2}^{l}$'s term is
\vspace{0.1cm}
\[
\lambda^{k+l} a_{kl}^{2}(\varphi(x)) = \lambda a_{kl}^{2}(x).
\vspace{0.1cm}
\]
By Proposition 3.8\,\,\,2), $a_{kl}^{2}(x) = 0$. Thus, by induction on $l$, $a_{kl}^{2}$ is identically equal to $0$ for all $k + l \geq 2$. Therefore, the function $\xi_{2}$ satisfying (L2) takes the form
\vspace{0.1cm} 
\[
\xi_{2}(z_{1},z_{2},x) = z_{2} + a_{00}^{2}(x), \,\,\,  a_{00}^{2} \in \mathcal{Z}_{\varphi,\lambda}.
\vspace{0.1cm}
\]

Next, we determine the function $\xi_{1}$ which satisfies (L1). The functional equation of (L1) is reduced to
\[
\sum_{k + l \geq 0} a_{kl}^{1}(\varphi(x)) (\lambda z_{1} + z_{2})^{k} (\lambda z_{2})^{l} = \sum_{k + l \geq 0} \lambda a_{kl}^{1}(x) z_{1}^{k} z_{2}^{l} + (z_{2} + a_{00}^{2}(x)).
\]
On $z_{1}^{0}z_{2}^{0}$'s term ($k = l = 0$), the equation is 
\vspace{0.1cm}
\[
a_{00}^{1}(\varphi(x)) = \lambda a_{00}^{1}(x) + a_{00}^{2}(x).
\vspace{0.1cm}
\]
On $z_{1}^{1}z_{2}^{0}$'s term ($k = 1$, $l =0$), the equation is 
\vspace{0.1cm}
\[
\lambda a_{10}^{1}(\varphi(x)) = \lambda a_{01}^{1}(x).
\vspace{0.1cm}
\]
By Proposition 3.8\,\,\,1), $a_{10}^{1}(x) = a_{10}^{1}(0) = 1$. On $z_{1}^{0}z_{2}^{1}$'s term ($k = 0$, $l = 1$), the equation is
\vspace{0.1cm}
\[
a_{10}^{1}(x) + \lambda a_{01}^{1}(\varphi(x)) = \lambda a_{01}^{1}(x) + 1.
\vspace{0.1cm}
\]
Since $a_{10}^{1}$ is equal to $1$, above equation is
\vspace{0.1cm}
\[
a_{01}^{1}(\varphi(x)) = a_{01}^{1}(x).
\vspace{0.1cm}
\]
By Proposition 3.8\,\,\,1), $a_{01}^{1}(x) = a_{01}^{1}(0) = 0$. If $k + l \geq 2$, the equation of $z_{1}^{k}z_{2}^{l}$'s term is
\[
\sum_{j=0}^{l}
\left(
\begin{array}{c}
k+l-j
\\
k
\end{array}
\right)
\lambda^{k+j} a_{k+l-j\,j}^{1}(\varphi(x)) = \lambda a_{kl}^{1}(x),
\]
which is same as $a_{kl}^{2}$ for $k+l \geq 2$, so that $a_{kl}^{1}$ is identically equal to 0 for all $k + l \geq 2$. Therefore, $g$ belongs to $\mathit{Aut}\,(\tilde{R},\tilde{H};T)$ iff it takes the following form.
\vspace{0.1cm}
\[
g (z_{1},z_{2},x) = \left( z_{1} + a_{00}^{1}(x),\, z_{2} + a_{00}^{2}(x),\, \eta(x) \right),
\vspace{0.1cm}
\]
where $a_{00}^{1}$, $a_{00}^{2} \in C^{\infty}([0,\infty);\mathbb{C})$  satisfies the functional equations
\vspace{0.1cm}
\[
a_{00}^{1}(\varphi(x)) = \lambda a_{00}^{1}(x) + a_{00}^{2}(x), \,\,\, a_{00}^{2}(\varphi(x)) = \lambda a_{00}^{2}(x)
\vspace{0.1cm}
\]
\textit{i.e.} $(a_{00}^{1},a_{00}^{2}) \in \mathcal{S}_{\varphi,\lambda}$ and $\eta \in Z_{\varphi}$.
\vspace{0.1cm}

Consequently, we obtain the following.
\\
\\
\textbf{Theorem 3.12} \,\, If $\tau = 1$ and $\lambda = \mu$, then the group $\textit{Aut}(R,H)$ admits a following sequence of extension,
\vspace{0.1cm}
\[
0 \,\to\, \mathcal{S}_{\varphi,\lambda} \,\to\, \mathit{Aut}\,(R,H) \cong \mathcal{S_{\varphi,\lambda}} \rtimes Z_{\varphi} \,\to\, Z_{\varphi} \,\to\, 1 
\vspace{0.1cm}
\]
where $\mathcal{S}_{\varphi,\lambda}$ is an infinite dimensional vector space described in Section 2 which is the set of functions $(a_{00}^{1},a_{00}^{2})$ in Case 4 and $\eta \in Z_{\varphi}$ acts on $(\beta_{1},\beta_{2}) \in \mathcal{S}_{\varphi,\lambda}$ by $(\beta_{1}(x),\beta_{2}(x)) \mapsto (\beta_{1}(\eta(x)),\beta_{2}(\eta(x))$. 
\\

\quad \,  \textbf{Case 5}. \, $\tau \not= 0$ \,and\, $\lambda = \mu^{p}$ for some $p \geq 2$.
\\
\\
\quad \, Let $S$ be the diffeomorphism of $\tilde{R}$ given by $S(z_{1},z_{2},x) = (z_{1}, \tau^{\frac{1}{p}} z_{2},x)$. Then, for any $(z_{1},z_{2},x) \in \tilde{R}$,
\vspace{0.1cm}
\[
S \circ T \circ S^{-1} (z_{1},z_{2},x) = (\mu^{p} z_{1} + z_{2}^{p}, \mu z_{2}, \varphi(x)).
\vspace{0.1cm}
\]
Therefore, a LC Reeb component $R = \tilde{R} / T^{\mathbb{Z}}$ is leafwise holomorphic foliated diffeomorphic to a LC Reeb component $R^{\prime} = \tilde{R} / (S \circ T \circ S^{-1} )^{\mathbb{Z}}$, so we may assume that $\tau$ is equal to 1.

If $\tau = 1$ and ${\lambda = \mu^{p}}$, Condition (L) is written as 
\[
(\textrm{L1}) \quad \xi_{1}(\mu^{p} z_{1} + z_{2}^{p},\, \mu z_{2},\, \varphi(x)) = \mu^{p} \xi_{1}(z_{1},z_{2},x) + \xi_{2}(z_{1},z_{2},x)^{p}, 
\vspace{0.2cm}
\]
\[
(\textrm{L2}) \quad \xi_{2}(\mu^{p} z_{1} + z_{2}^{p},\, \mu z_{2},\, \varphi(x)) = \mu \xi_{2} (z_{1},z_{2},x).
\vspace*{0.3cm}
\]
First, we determine the function $\xi_{2}$ which satisfies (L2). The functional equation of (L2) is reduced to
\[
\sum_{k + l \geq 0} a_{kl}^{2}(\varphi(x)) (\mu^{p} z_{1} +  z_{2}^{p})^{k} (\mu z_{2})^{l} = \sum_{k + l \geq 0} \mu a_{kl}^{2}(x) z_{1}^{k} z_{2}^{l}.
\]
Furthermore, we compare both coefficient functions of $z_{1}^{k}z_{2}^{l}$'s term.

Let $p$ be the quotient and $r$ be the remainder when $l$ is divided by $p$, \textit{i.e.} $l = pq + r$, $q \geq 0$, $0 \leq r \leq q-1$. If $k \geq 1$ and $l \geq 0$, then the functional equation of $z_{1}^{k}z_{2}^{l}$' term is 
\[
\sum_{j=0}^{q}
\left(
\begin{array}{c}
k+j
\\
j
\end{array}
\right)
\mu^{pk} \mu^{p(q-j)+r} a_{k+j\,p(q-j) + r}^{2}(\varphi(x)) = \mu a_{k\,pq+r}^{2}(x).
\]
If $k \geq 1$ and $q=0$ ($l = r$), then above equation is
\vspace{0.1cm}
\[
\mu^{pk} \mu^{r} a_{kr}^{2}(\varphi(x)) = \mu a_{kr}^{2}(x).
\vspace{0.1cm}
\]
By Proposition 3.8\,\,\,2), $a_{kr}^{2}(x) = 0$ for $k \geq 1$ and $0 \leq r \leq p-1$. For $k \geq 1$ and $q \geq 1$, if $a_{k + q\,r}^{2},\,a_{k + q -1 \, q + r}^{2},\, \cdots \,$ and $a_{k+1\,p(q-1)+r}^{2}$ are identically equal to $0$, then the equation of $z_{1}^{k}z_{2}^{pq+r}$'s term is
\vspace{0.1cm}
\[
\mu^{pk} \mu^{pq + r} a_{k\,pq + r}^{2}(\varphi(x)) = \mu a_{k\,pq+r}^{2}(x).
\vspace{0.1cm}
\]
By Proposition 3.8\,\,\,2), $a_{k\,pq+r}^{2}(x) = 0$. Thus, by induction on $q$, $a_{kl}^{2}$ is identically equal to $0$ for $k \geq 1$ and $l \geq 0$. If $k=0$ and $l \geq 0$, then the equation of $z_{1}^{0}z_{2}^{l}$'s term is
\vspace{0.1cm}
\[
\sum_{j = 0}^{q} \mu^{p(q-j)+r} a_{j\,p(q-j) + r}^{2}(\varphi(x)) = \mu a_{0\,pq + r}^{2}(x). 
\vspace{0.1cm}
\]
If $q = r = 0$, then above equation is
\vspace{0.1cm}
\[
a_{00}^{2} (\varphi(x)) = \mu a_{00}^{2}(x), \,\,\textit{i.e.}\,\, a_{00}^{2} \in \mathcal{Z}_{\varphi,\mu}.
\vspace{0.1cm}
\]
If $q = 0$ and $r \not = 0$, then the equation of $z_{1}^{0}z_{2}^{r}$'s term is
\vspace{0.1cm}
\[
\mu^{r} a_{0r}^{2}(\varphi(x)) = \mu a_{0r}^{2}(x).
\vspace{0.1cm}
\]
If $r = 1$, then $a_{01}^{2}(x) = a_{01}^{2}(0) = 1$ by Proposition 3.8\,\,\,1), and if $2 \leq r \leq p-1$, then $a_{0r}^{2}(x) = 0$ by Proposition 3.8\,\,\,2). If $q \geq 1$, $a_{1\,p(q-1) + r}^{2},\,\cdots,\,a_{q-1\,p+r}^{2}$ and $a_{qr}^{2}$ are identically equal to 0, and so the equation of $z_{1}^{0}z_{2}^{pq+r}$'s term is
\vspace{0.1cm}
\[
\mu^{pq + r} a_{0\,pq +r}^{2}(\varphi(x)) = \mu a_{0\,pq + r}^{2}(x).
\vspace{0.1cm}
\]
By Proposition 3.8\,\,\,2), $a_{0\,pq+r}^{2}(x) = 0$. Thus, $a_{0l}^{2}$ is identically equal to 0 for $l\geq 0$. Therefore, the function $\xi_{2}$ satisfying (L2) takes the form
\vspace{0.1cm} 
\[
\xi_{2}(z_{1},z_{2},x) = z_{2} + a_{00}^{2}(x), \,\,\,  a_{00}^{2} \in \mathcal{Z}_{\varphi,\mu}.
\vspace{0.1cm}
\]

Next, we determine the function $\xi_{1}$ satisfying (L1). The functional equation of (L1) is reduced to
\vspace{0.1cm}
\[
\sum_{k + l \geq 0} a_{kl}^{1}(\varphi(x)) (\mu^{p} z_{1} + z_{2}^{p} )^{k} (\mu z_{2})^{l} = \sum_{k + l \geq 0} \mu^{p} a_{kl}^{1}(x) z_{1}^{k} z_{2}^{l} + (z_{2} + a_{00}^{2}(x))^{p}.
\vspace{0.1cm}
\]
If $k \geq 1$ and $l \geq 0$, then the functional equation of $z_{1}^{k}z_{2}^{l}$' term is 
\[
\sum_{j=0}^{q}
\left(
\begin{array}{c}
k+j
\\
j
\end{array}
\right)
\mu^{pk} \mu^{p(q-j)+r} a_{k+j\,p(q-j) + r}^{1}(\varphi(x)) = \mu^{p} a_{k\,pq+r}^{1}(x).
\]
If $k \geq 1$ and $q=0$ ($l = r$), then above equation is
\vspace{0.1cm}
\[
\mu^{pk} \mu^{r} a_{kr}^{1}(\varphi(x)) = \mu^{p} a_{kr}^{1}(x).
\vspace{0.1cm}
\]
Then, if $(k,r) = (1,0)$, $a_{10}^{1}(x) = a_{10}^{1}(0) = 1$ by Proposition 3.8\,\,\,1), and if $(k,r) \not = (1,0)$, $a_{kr}^{1}(x) = 0$ by Proposition 3.8\,\,\,2). For $k \geq 1$ and $q \geq 1$, if $a_{k + q\,r}^{1},\,\cdots,\,a_{k + 2\,p(q-2)+r}$ and $a_{k+1\,p(q-1)+r}^{1}$ are identically equal to $0$, then the equation of $z_{1}^{k}z_{2}^{pq+r}$'s term is
\vspace{0.1cm}
\[
\mu^{pk} \mu^{pq + r} a_{k\,pq + r}^{1}(\varphi(x)) = \mu^{p} a_{k\,pq+r}^{1}(x).
\vspace{0.1cm}
\]
By Proposition 3.8\,\,\,2), $a_{k\,pq+r}^{1}(x) = 0$. Thus, by induction on $q$, $a_{kl}^{1}$ is identically equal to $0$ for $k \geq 1$, $l \geq 0$ and $k + l \geq 2$. If $k=0$ and $l \geq 0$, then the functional equation of $z_{1}^{0}z_{2}^{l}$'s term is the following.
\vspace{0.1cm}
\[
\mu^{l} a_{0l}^{1}(\varphi(x)) = \mu^{p} a_{0l}^{1}(x) + \left( \begin{array}{c} p \\ l \end{array} \right) a_{00}^{2}(x)^{p-l} \quad (\textrm{if}\,\,\, 0 \leq l \leq p-1)
\vspace{0.2cm}
\]
\[
a_{10}^{1}(\varphi(x)) + {\mu}^{p} a_{0p}^{1}(\varphi(x)) = \mu^{p} a_{0p}^{1}(x) + 1 \quad (\textrm{if}\,\,\,l = p)
\vspace{0.1cm}
\]
\[
\sum_{j = 0}^{q} \mu^{p(q-j)+r} a_{j\,p(q-j) + r}^{1}(\varphi(x)) = \mu^{p} a_{0\,pq + r}^{2}(x) \quad (\textrm{if}\,\,\,l \geq p + 1) 
\vspace{0.2cm}
\]
If\, $0 \leq l \leq p-1$, by dividing both sides by $\mu^{l}$, the equation is reduced to
\vspace{0.1cm}
\[
a_{0l}^{1} (\varphi(x)) = \mu^{p-l} a_{00}^{1}(x) + c_{l} a_{00}^{2}(x)^{p-l}, \,\,\,\, c_{l} = \left( \begin{array}{c} p \\ l \end{array} \right) \mu^{-l}.
\vspace{0.1cm}
\]
If $l = p$, the equation is reduced to
\vspace{0.1cm}
\[
a_{0p}^{1}(\varphi(x)) = a_{0p}^{1}(x),
\vspace{0.1cm}
\]
because $a_{10}^{1}$ is equal to 1. By Proposition 3.8\,\,\,1), $a_{0p}^{1}(x) = a_{0p}^{1}(0) = 0$. If $l \geq p+1$, $q$ is not equal to $0$. Then, $a_{1\,p(q-1)+r}^{1},\, \cdots,\, a_{q-1\,p+r}^{1}$ and $a_{qr}^{1}$ are identically equal to $0$, and so the equation is reduced to
\vspace{0.1cm}
\[
\mu^{pq+r} a_{0\,pq+r}^{1}(\varphi(x)) = \mu^{p} a_{0\,pq+r}^{1}(x).
\vspace{0.1cm}
\]
By Proposition 3.8\,\,\,2), $a_{0\,pq + r}^{1}(x) = 0$ for $l = pq + r \geq p+1$. Therefore, $g$ belongs to $\mathit{Aut}\,(\tilde{R},\tilde{H};T)$ iff it takes the following form.
\vspace{0.1cm}
\[
g (z_{1},z_{2},x) = \left( z_{1} + \sum_{j = 0}^{p-1} a_{0j}^{1}(x) z_{2}^{j},\, z_{2} + a_{00}^{2}(x),\, \eta(x) \right),
\vspace{0.1cm}
\]
where $\eta \in Z_{\varphi}$ and $a_{0j}^{1}\,\,(j = 0,\,1,\,\cdots\,,\,p-1)$, $a_{00}^{2} \in C^{\infty}([0,\infty);\mathbb{C})$ satisfy the functional equations
\vspace{0.1cm}
\[
a_{0j}^{1}(\varphi(x)) = \mu^{p-j} a_{0j}^{1}(x) + c_{j} a_{00}^{2}(x)^{p-j}, \,\,\, a_{00}^{2}(\varphi(x)) = \mu a_{00}^{2}(x),
\vspace{0.2cm}
\]
where $c_{j} = \left( \begin{array}{c} p \\ j \end{array} \right) \mu^{-j} \in \mathbb{C}$.
\vspace{0.4cm}
\\

Let us consider the following system of functional equations on $\beta_{1,0},\,\cdots\,,\,\beta_{1,p-1}$\, and $\beta_{2}$ $\in C^{^{\infty}}([0,\infty);\mathbb{C})$ concerning $\varphi$ and $\mu$.
\vspace{7pt}
\\
\quad Equation (I\hspace{-.1em}I\hspace{-.1em}I) :
\quad $\beta_{1,0}(\varphi(x)) = \mu^{p} \beta_{1,0}(x) + c_{0} \beta_{2}(x)^{p}$,
\vspace{5pt}
\\
\hspace{101pt}
$\beta_{1,1}(\varphi(x)) = \mu^{p-1} \beta_{1,1}(x) + c_{1} \beta_{2}(x)^{p-1}$,
\vspace{2pt}
\\
\hspace{120pt}
\vdots
\vspace{2pt}
\\
\hspace{101pt}
$\beta_{1,p-1}(\varphi(x)) = \mu \beta_{1,p-1}(x) + c_{p-1} \beta_{2}(x)$,
\vspace{5pt}
\\
\hspace{101pt}
$\beta_{2}(\varphi(x)) = \mu \beta_{2}(x)$. 
\vspace{7pt}
\\
where $c_{0},\,\cdots,\,c_{p-1} \in \mathbb{C}$ are non-zero constants. Let us take the space $\mathcal{S}_{\varphi,\mu}(\mathbf{c})$ of solutions to Equation (I\hspace{-.1em}I\hspace{-.1em}I) on $(0,\infty)$, where $\mathbf{c} = (c_{0},\,\cdots,\,c_{p-1}) \in \mathbb{C}^{p}$. If we assign $\beta_{2}$ to a solution $(\beta_{1,0},\,\cdots,\,\beta_{1,p-1},\,\beta_{2}) \in \mathcal{S}_{\varphi,\mu}(\mathbf{c})$, we obtain the projection $P_{2} : \mathcal{S}_{\varphi,\mu}(\mathbf{c}) \to \mathcal{Z}_{\varphi,\mu}$. Here the kernel of $P_{2}$ is nothing but $\mathcal{Z}_{\varphi,\mu^{p}} \times \cdots \times \mathcal{Z}_{\varphi,\mu}$. If we fix any solution $\beta^{*} \in \mathcal{Z}_{\varphi,\mu}$ which never vanishes, we also see that the projection $P_{2}$ is surjective because for any $\beta_{2} \in \mathcal{Z}_{\varphi,\mu}$
\vspace{0.1cm}
$$
\beta_{1,j}(x)= \frac{c_{j}}{\mu^{p-j} \log \mu^{p-j}} \beta_2(x)^{p-j} \log \{\beta^*(x)^{p-j}\} \,\,\, (j = 0,\,1,\,\cdots,\,p-1)
\vspace*{0.1cm}
$$
gives a solution $(\beta_{1,0},\,\cdots,\,\beta_{1,p-1},\,\beta_{2}) \in \mathcal{S}_{\varphi,\mu}(\mathbf{c})$, where for $\log \{\beta^{*}(x)^{p-j}\}$ for any smooth branch can be taken. Therefore, as a vector space, $\mathcal{S}_{\varphi,\mu}(\mathbf{c})$ has a structure that
\vspace{0.1cm}
$$
0 \,\to\, \mathcal{Z}_{\varphi,\mu^{p}} \times \cdots \times \mathcal{Z}_{\varphi,\mu} \,\to\, \mathcal{S}_{\varphi,\mu}(\mathbf{c}) \,\to\, \mathcal{Z}_{\varphi,\mu} \,\to\, 0
\vspace{0.1cm}
$$
is a short exact sequence.

Any solution $(\beta_{1,0},\,\cdots,\,\beta_{1,p-1},\,\beta_{2}) \in \mathcal{S}_{\varphi,\mu}(\mathbf{c})$ extends to $[0,\infty)$ so as to be  smooth functions which are flat at $x=0$ because $(\beta_{1,j},\,c_{j}\beta_{2}^{p-j}) \in \mathcal{S}_{\varphi,\mu^{p-j}}$ for $j = 0,\,1\,\cdots,\,p-1$ (see Theorem 3.1 in [HM]). 
\\
\\
\textbf{Theorem 3.13} \,\, If $\tau = 1$ and $\lambda = \mu^{p}$ for some $p \geq 2$, then the group $\textit{Aut}(R,H)$ admits a following sequence of extension,
\vspace{0.1cm}
\[
0 \,\to\, \mathcal{S}_{\varphi,\mu}(\mathbf{c}) \,\to\, \mathit{Aut}\,(R,H) \,\to\, Z_{\varphi} \,\to\, 1 
\vspace{0.1cm}
\]
where $\mathcal{S}_{\varphi,\mu}(\mathbf{c})$ is an infinite dimensional vector space which is the set of functions $(a_{00}^{1},\,\cdots,\,a_{0p-1}^{1},\,a_{00}^{2})$ in Case 5 and $\mathbf{c} = (c_{0},\,\cdots,\,c_{p-1}) \in \mathbb{C}^{p}$, $c_{j} = \left( \begin{array}{c} p \\ j \end{array} \right) \mu^{-j}$.
\vspace{0.2cm}

\section{\large{Higher dimensional Reeb components}}

To close this article, we make some remarks on the automorphisms of higher dimensional LC Reeb components. Let $\tilde{R}$ be $\mathbb{C}^{n} \times [0,\infty) \backslash \{ (O,0) \}$ and $\varphi \in \textit{Diff}^{\,\,\infty}([0,\infty))$ be a diffeomorphism which is tangent to the identity to the infinite order at $0$ and satisfies $\varphi(x) - x > 0$ for $x > 0$. Also take an automorphism $G \in \textit{Aut}\,(\mathbb{C}^{n},\{O\})$ which is expanding. Let $T$ be a diffeomorphism given by $T = G \times \varphi$. Then we obtain a LC Reeb component $(R,\mathcal{F},J) = (\tilde{R}, \tilde{\mathcal{F}},J_{\textrm{std}}) / T^{\mathbb{Z}}$ as the quotient, as well as the boundary Hopf manifold $H = \mathbb{C}^{n} \backslash \{ O \} / G^{\mathbb{Z}}$. 

Let $\textit{Aut}\,R$ be the group of leafwise holomorphic smooth automorphisms of $R$ and $\textit{Aut}\,H$ be the group of holomorphic automorphisms of the boundary Hopf manifold $H$. Then, by the similar argument, we can obtain the claims as Proposition 3.2, 3.3 and 3.7.
\\
\\
\textbf{Proposition 4.1} \, For $n \geq 1$, $\textit{Aut}\,R$ and $\textit{Aut}\,(R,H)$ are isomorphic to $\textit{Aut}\,(\tilde{R};T)/T^{\mathbb{Z}}$ and $\textit{Aut}\,(\tilde{R},\tilde{H};T)$ respectively. 
\\
\\
\textbf{Proposition 4.2} \, For $n \geq 2$, the restriction map $r_{H}:\textit{Aut}\,R \to \textit{Aut}\,H$ is surjective. 
\\

If $G$ is in some normal forms like in Theorem 1.2, then we can most probably compute the group $\textit{Aut}\,H$ and the kernel $\textit{Aut}\,(R,H)$ of the restriction map $r_{H}$ for $n > 2$. However, in order to determine these groups in all cases, we need to classify Hopf manifolds and to give normal forms like Theorem 1.2 for $n = 2$.  
\\
\\
\textbf{Example 4.3} \,\, Let us look at the simplest case, where $G$ is a diagonal matrix $\lambda I_{n} \in GL(n,\mathbb{C})$ with $|\lambda| > 1$. Then the automorphism group $\textit{Aut}\,R$ admits the following sequence of extensions.
\vspace{0.1cm}
\[
0 \,\to\, \textit{Aut}\,(R,H) \,\to\, \textit{Aut}\,R \,\to\, \textit{Aut}\,H \cong GL(n,\mathbb{C}) / \{ \lambda I_{n}\}^{\mathbb{Z}} \,\to\, 0,
\vspace{0.1cm}
\] 
\[
0 \,\to\, (\mathcal{Z_{\varphi,\lambda}})^{n} \,\to\, \textit{Aut}\,(R,H) \cong (\mathcal{Z}_{\varphi,\lambda})^{n} \rtimes Z_{\varphi}  \,\to\, Z_{\varphi} \,\to\, 1
\vspace{0.3cm}
\]
where $\eta \in Z_{\varphi}$ acts on $(\beta_{1},\,\cdots,\,\beta_{n}) \in (\mathcal{Z}_{\varphi,\lambda})^{n}$ by $(\beta_{1}(x),\,\cdots,\,\beta_{n}(x)) \mapsto (\beta_{1}(\eta(x)),\,\cdots,\,\beta_{n}(\eta(x)))$.

\vspace{0.5cm}
\begin{flushright}
Tomohiro HORIUCHI
\\
\mbox{}
\\
{\small Department of Mathematics, Chuo University
\\
1-13-27 Kasuga Bunkyo-ku, Tokyo, 112-8551, Japan
\mbox{}
\\
e-mail : horiuchi@gug.math.chuo-u.ac.jp
\\
}
\end{flushright}

\end{document}